\UseRawInputEncoding
\documentclass[11pt]{article}

\usepackage{multicol}
\usepackage{amsmath,amssymb,amsfonts,amsthm,mathrsfs,stmaryrd}
\usepackage[usenames,dvipsnames]{color}

\newcommand{\textcite}{\cite}
\newcommand{\autocite}{\cite}

\theoremstyle{plain}
\newtheorem{thm}{Theorem}[section]
\newtheorem{lem}[thm]{Lemma}
\newtheorem{prop}[thm]{Proposition}

\theoremstyle{definition}

\theoremstyle{remark}

\newcommand\prf{\noindent{\bf Proof}}

\newcommand{\ws}{\overset{\operatorname{w}_1}{\longrightarrow}}

\title{\small\bf CRITICAL FLUID LIMIT OF A GATED\\  PROCESSOR SHARING QUEUE}

\author{ {\small\sc H.\ Christian Gromoll, \ Katelynn D.\ Kochalski} \\
{\em\footnotesize University of Virginia and SUNY Geneseo} }

\begin{document}
\maketitle

\begin{abstract}
We consider a sequence of single-server queueing models operating under a
service policy that incorporates batches into processor sharing: arriving jobs
build up behind a gate while waiting to begin service, while jobs in front of
the gate are served according to processor sharing. When they have been
completed, the waiting jobs move in front of the gate and the cycle repeats. 
We model this system with a pair of measure valued processes describing the
jobs in front of and behind the gate. Under mild asymptotically critical 
conditions and a law-of-large-numbers scaling, we prove that the
pair of measure-valued processes converges in distribution to an easily
described limit, which has an interesting periodic dynamics.
\end{abstract}



\section{Introduction}

Consider a single server queue operating under a service policy that is a
variant of processor sharing: instead of all jobs receiving simultaneous,
equal service regardless of arrival time, jobs are batched so that only
members of the active batch receive simultaneous equal service, while newly
arriving jobs wait and receive no service until the current batch has
completed processing; the next batch consists of all jobs that have arrived
and are waiting at the completion time of the previous batch. One can think of
this policy as a gated form of processor sharing, in which the gate superimposes
some low-resolution aspects of the First-In-First-Out (FIFO) policy onto pure
processor sharing. 

One aspect of processor sharing that has been criticized is that it can take a
long time for large jobs to complete service if there are many small jobs
present that slow down the system.  Queues under FIFO don't encounter this
problem.  One way to mitigate such slowdown while retaining essential features
of processor sharing is to incorporate a gate that divides the arrival stream
into batches that are sequenced according to FIFO, but processed individually
via processor sharing.

Gated service policies have been studied in many different settings.  The
works \autocite{RegeSengupta}, \autocite{Avi-ItzhakHalfin}, and
\autocite{JagermanSengupta} consider gated processor sharing systems allowing at
most $m<\infty$ jobs per batch.  Rege and Sengupta \textcite{RegeSengupta} analyze the
distributions of queue length, mean sojourn time, and busy periods for Poisson
arrivals and exponential service times.  Avi-Itzhak and Halfin \textcite{Avi-ItzhakHalfin} study the
sojourn and response time distributions for general service times, and
Jagerman and Sengupta \textcite{JagermanSengupta} allow bulk Poisson arrivals, general service
times, and study the waiting time, queue length, and batch size distributions
for this setting. 

Gated service policies have also been studied in polling models (systems with
a single server and multiple queues being processed cyclically by the server).
Boxma, Kella, and Kosinski \textcite{BoxmaKellaKosinski} study the distribution of queue length and
workload for a system with Poisson arrivals.  Fuhrmann, Levy, and Sidi \textcite{FuhrmannLevySidi}
derive first and second moments of  each queue length as well as the expected
customer wait time under Poisson arrivals. These models are related to the
system we consider here which could be viewed as an idealization of a gated
polling model. 

A rather different service policy combining FIFO and processor sharing
protocols is head-of-the-line proportional processor sharing. In this model,
there are jobs of several different classes and the server uses processor
sharing to simultaneously serve one job from each class. Within each class,
the jobs are served according to FIFO. Fluid and diffusion limits under heavy
traffic conditions were established for this model by Bramson \textcite{Bramson} and
Williams \textcite{WilliamsR}.

Another policy that compensates for system slowdown is referred to as limited
processor sharing. In this model, processor sharing is used to serve maximally
$k$ jobs. If there are less than $k$ jobs in service, an arriving job begins
service immediately. If there are $k$ jobs in service, arriving jobs must wait
to begin service. Once a job departs, the job that has been waiting the
longest begins service. Dai, Zhang, and Zwart\textcite{DaiZhangZwart} derive fluid limits for such
a system with general arrival and service time distributions.

The main difference between the model studied here and most of the systems
described above is that there is no limit on the number of jobs that could be
receiving simultaneous service. The gate simply remains closed until the batch
in service is complete and then opens to let all waiting jobs begin service.
Consequently, the sizes of batches are affected by the processing time of the
previous batch, which determines the length of time for arriving jobs to build
up the next batch. This results in an interesting feedback mechanism that
creates dependency far out into the future. Understanding the long-run
behavior of this model is therefore non-trivial.  

We note that another name for the service policy considered here could be
``Most Attained Service,'' contrasting it with the Least Attained Service
policy. Indeed if the server only
works on jobs with the most attained service, employing processor sharing to
simultaneously process ties, one has the same protocol as already described.

While the total system workload is the same for this system as for any other
single-server queue under a work-conserving policy, the behaviors of other
important performance measures such as queue length or sojourn time are not
obvious. The goal of this article is to propose a stochastic model of the gated
processor sharing queue that is detailed enough to shed light on such
performance measures via approximating fluid models. Building on the
well-established approach for tracking processor-sharing-style policies, we
use measure-valued state descriptors to track all residual service times and
we formulate the limiting fluid model as measure-valued dynamics. Under mild
assumptions allowing general interarrival and service times, we show that the
stochastic model converges to the fluid model under a law of large numbers
scaling. 

Apart from the performance conclusions that can be gleaned from our fluid
limit, a principal motivation for studying this model is to explore
methods for handling systems with cyclic fluid limits, particularly as
pertains to the relationship between fluid and diffusion limits. Diffusion
limits are second-order approximations to queueing systems that provide richer
descriptions of their inherent randomness than do first-order fluid limits.
The two types of result are intimately related however. A well-established
technique for proving diffusion limits is to combine multiple overlapping
sections of fluid limits with a steady state result about them, known as state
space collapse. Namely, if the fluid limits converge (as time tends to
infinity) to certain invariant states, then the overlapping fluid limits often
combine to give a limiting diffusion process on the set of invariant states;
see for example \cite{Grom2004}. 

As will be shown in this article, the fluid limits we establish do not
converge asymptotically in time to fixed invariant states. In fact, except for
the zero state, our model has no invariant states. This calls into question
the prospect of proving diffusion limits with the above approach. It turns out
however that the fluid limit of our model has certain fixed orbits, such that
starting the fluid model within an orbit results in perpetual cycling through
it. Moreover, each initial condition has a corresponding orbit to which the
fluid model will converge, in the sense of getting asymptotically close to an
appropriately orbiting state within it.

This behavior resurrects the possibility of establishing diffusion limits via
fluid limits, albeit in a generalized form. A limiting diffusion process would
have to be defined on the set of orbits and combined with some kind of
averaging principle within orbits, in order properly approximate the prelimit
system. This will be the subject of future work. 

To our knowledge, the literature has not given much attention to
asymptotically orbiting fluid limits and their relation to potentially
generalized diffusion limits. This seems to be interesting ground for further
research.  Motivated by the above, the goal of the present article is to
define a dynamic model that captures the distinctiveness of gated processor
sharing, and to prove that under a law of large numbers scaling, this model
converges to a descriptive fluid model.

\subsection{Notation}

Let $\mathbb{Z}_{+}=\{0, 1, 2,\ldots\}$ and $\mathbb{R}_+=[0,\infty)$. Denote
the integer part of $a$ by $\lfloor a\rfloor$ and the negative part of
$a\in\mathbb{R}$ by $a^-=-a\vee 0$. For $w>0$, let $\llbracket
t\rrbracket_w=t-\lfloor t/w\rfloor w$ denote the residue modulo $w$ of
$t\in\mathbb{R}_+$, and let $\llbracket t\rrbracket_0=0$.  The identity
function on $\mathbb{R}_+$ is denoted $\chi$ and write $\chi_a=\chi\wedge a$
and $\chi_{[a,\infty)}=\chi 1_{[a,\infty)}$. We use the convention that
sums of the form $\sum_{i=a}^{b}s_{i}$ equal zero if $b<a$.

Let $\mathbf{C}=\mathbf{C}(\mathbb{R}_+)$ be the space of continuous
real-valued functions on $\mathbb{R}_+$ and $\mathbf{C}_b\subset\mathbf{C}$
the subspace of bounded continuous functions. Let $\mathcal{M}_1$ be the space
of finite non-negative Borel measures on $\mathbb{R}_+$ with finite first
moment. For an integrable $f:\mathbb{R}_+\to\mathbb{R}$ and
$\zeta\in\mathcal{M}_1$ we write \[ \langle f,\zeta\rangle =
\int_{\mathbb{R}_+}fd\zeta. \] We will frequently ``shift $\zeta$ left by $x$
and remove any mass at or below zero.'' The resulting measure is is written
$\zeta(\cdot+_0x)$ and means \[ \langle f,\zeta(\cdot+_0x)\rangle = \langle
(1_{(0,\infty)}f)(\cdot-x),\zeta\rangle, \] where it is understood that $f$ is
always extended to equal zero on $(-\infty,0)$ so that the above is
well-defined on $\mathbb{R}_+$. 

We give $\mathcal{M}_1$ the structure of a Polish space by endowing it with
the Wasserstein$_1$ topology: $\zeta_n\ws\zeta$ if and only if
$\langle\chi,\zeta_n\rangle \to \langle \chi,\zeta\rangle$ and $\langle
f,\zeta_n\rangle \to \langle f,\zeta\rangle$ for all $f\in\mathbf{C}_b$. That
is, $\operatorname{w}_1$-convergence (also called Kantorovich convergence) is
equivalent to weak convergence plus convergence of first moments. The space
of paths in $\mathcal{M}_1$ or $\mathbb{R}_+$ that are right-continuous with
left limits is denoted $\mathbf{D}([0,\infty),\mathcal{M}_1)$ or
$\mathbf{D}([0,\infty),\mathbb{R}_+)$ respectively. These spaces have the
usual Skorohod $J_1$-topology.

We write $X \sim Y$ if $X$ and $Y$ have the same distribution and denote
convergence in distribution by  $\Rightarrow$.

\section{Model and main result}

We begin by defining the stochastic model for a gated processor sharing queue,
defining the analogous fluid model, and stating the main result establishing
the fluid model as the limit of a sequence of scaled stochastic models.

\subsection{Stochastic model}\label{s.stochasticModel}

Exogenous arrivals to the system are given by a (possibly delayed) renewal
process $E(\cdot)$ with finite mean arrival rate $\alpha$ and $E(0)=0$.
$E(t)$ is the number of jobs that have arrived to the system by time $t\ge0$, not
counting any inital jobs that may be present at  $t=0$.



The service times of arriving jobs are taken from a sequence $\{v_{i}\}$ of
independent, identically distributed random variables with distribution $\nu$,
a Borel probability measure on $\mathbb{R}_{+}$ with $\nu(\{0\})=0$ and finite
mean $1/\beta$.

For the initial condition, let $Z_{0}$ be a non-negative integer valued random
variable, and let $\{\tilde{v}_{j}\}$ be a sequence of positive random
variables. $Z_{0}$ is the initial queue length and $\{\tilde{v}_{j} : j=1, 2,
\ldots, Z_{0}\}$ are the service times of these initial jobs. Assume that 
$\mathbb{E}\left[\sum_{j=1}^{Z_{0}}\tilde{v}_{j}\right]<\infty$. 

The arrival process, service times, and initial condition above constitute the
stochastic primitives of the model. Note that service times and interarrival
times are generally distributed, so this is a G/G/1 type queuing model. Note
also that the general definition of the initial condition combined with the
allowance of a delayed renewal process, enables the model to describe systems
that, at time zero, may have already been operating in the past, and may thus
for example have initial service times and interarrival time distributed as
residuals of the nominal ones. 

It will be convenient to encode the primitives using elements of
$\mathcal{M}_1$.  Let
\[ \mathcal{B}_0=\sum_{j=1}^{Z_0}\delta_{\tilde v_j}, \qquad 
\mathcal{E}(t) = \sum_{i=1}^{E(t)}\delta_{v_i}, \quad t\ge 0, \]
where $\delta_v$ is the standard Dirac mass at $v$.

We next define several performance processes and other data that arise from
the stochastic primitives. The workload process is defined
\begin{equation}\nonumber
  W(t)=W_{0}+\langle\chi,\mathcal{E}(t)\rangle-t+I(t), \qquad t\ge0,
\end{equation}
where $W_{0}=\langle\chi,\mathcal{B}_0\rangle$
is the initial workload and 
\[ I(t)=\sup_{s\in[0,t]}\left(W_{0}+\langle\chi,\mathcal{E}(s)\rangle
-s\right)^-,\qquad t\ge 0, \]
is the cumulative idle time in $[0,t]$. Note that $W(t)$ is the same process
for all single-server queues under a work-conserving policy, and describes the
amount of time needed to empty the system if no further arrivals occurred. 

The distinguishing feature of gated processor sharing is that jobs are grouped
by arrival time into consecutive batches that can be defined from the workload
process. We now define the start time of each batch, as
well as finite Borel measures giving the profile of starting service times for
each batch. Let $\beta_{0}=0$ and for positive integers $k$ define $\beta_{k}$
inductively by
\begin{equation}\label{e.batchStartTimes}
 \beta_{k}=\inf\{s\geq \beta_{k-1}+W(\beta_{k-1}) \, : \, W(s)>0\}. 
\end{equation}
Note that almost surely, all $\beta_k$ are finite since $E(\cdot)$ is a
renewal process and by construction,  $\beta_k>\beta_{k-1}$ for all $k\ge1$.
For each $k\ge0$, $\beta_{k}$ gives the start time for the $k$th batch, where
batch zero refers to initial jobs whose service time profile is given by
$\mathcal{B}_0$, which may be zero if there are none. In this case the system
starts empty and the first idle period $[0,\beta_1)$ is associated with this
trivial batch zero. In all other cases, idle periods are associated with the
previously completed batch (that batch is still considered ``active'' even
  though it has completed processing), and are given by
  $[\beta_k+W(\beta_k),\beta_{k+1})$.  These intervals are clearly empty if
  there are jobs waiting at time $\beta_k+W(\beta_k)=\beta_{k+1}$, in which
  case there is no idle period between batches $k$ and $k+1$.

For $k\ge1$ define the $k$th batch profile by 
\begin{equation}
  \mathcal{B}_k=\mathcal{E}(\beta_k)-\mathcal{E}(\beta_{k-1}), \quad
k\ge 1. 
\end{equation}
Define the start time of the currently active batch
\begin{equation}
\beta(t)=\max\{\beta_k  :  \beta_k\le t\}, \qquad t\ge 0,
\end{equation}
the index of the currently active batch
\begin{equation}
\ell(t)=\max\{j  :  \beta_{j}\le t\}, \qquad t\ge0,
\end{equation}
and the starting profile of the currently active batch
\begin{equation}\label{e.batchProfiles}
\mathcal{B}(t)=\mathcal{B}_{\ell(t)}, \qquad t\ge0. 
\end{equation}
Note that $\beta(\beta(t))=\beta(t)$ and $\ell(\beta(t))=\ell(t)$ for all
$t\ge0$.

In gated processor sharing, the currently active batch is served according to
processor sharing, while any jobs that arrived during $(\beta(t),t]$ receive
no service until the next start time. The state of the system at time $t\ge0$
can therefore be split into two parts, a finite Borel measure $\sigma(t)$
describing the current profile of residual service times of all jobs in the
active batch, and a second finite Borel measure $\mu(t)$ describing the
service time profile of waiting jobs. 

Each job in the active batch receives simultaneous service at rate equal to
the inverse size of the batch $1/\langle 1,\sigma(t)\rangle,$ understood to be
zero if $\sigma(t)=0$ and the system is idling. That is, the measure
$\sigma(t)$ shifts to the left at varying rates, losing any mass that reaches
zero. Because newly arriving jobs do not slow this shifting down, its dynamics
can be described using a work conservation principle as follows.

Given $\zeta \in \mathcal{M}_1$, define $F_{\zeta}:\mathbb{R}_+\rightarrow
\mathbb{R}$ by
\[ F_\zeta(x)=  \langle \chi_x, \zeta\rangle, \qquad x\ge0, \] 
where $\chi_x=\chi \wedge x.$ If $\zeta=\sum_v\delta_v$ represents a profile
of service times, then $F_\zeta(x)$ represents the total amount of work the
processor sharing server must complete to provide every job with service time
$v$ an amount of service equal to $v\wedge x$, or equivalently, to shift
$\zeta$ left by $x$ and remove mass at or below zero. 

Consider the active batch at time $t\ge0$, which
was started at time $\beta(t)$ with starting service time profile
$\mathcal{B}(t)$. If $t\le \beta(t)+W(\beta(t))$, then it has been in
continuous service for time $t-\beta(t)$. The profile of residual service
times of the active batch will be given by $\mathcal{B}(t)$ shifted left by a
certain amount $S(t)$ with mass at or below zero removed. Since the server
hasn't idled, $S(t)$ must satisfy the relation
\[ t-\beta(t)=F_{\mathcal{B}(t)}(S(t)), \qquad t\in
  [\beta(t),\beta(t)+W(\beta(t)).\]

Let $x^*$ denote the supremum of the support of $\zeta$ and observe that
$F_{\zeta}$ is continuous and strictly increasing on $[0,x^*)$, and in the
case of bounded support $x^*<\infty$ is constant and equal to
$\langle\chi,\zeta\rangle$ on $[x^*,\infty).$ In case of unbounded support
$x^*=\infty$ and $F_{\zeta}$ approaches this constant asymptotically. In either
case, $F_{\zeta}$ has a well-defined continuous increasing inverse
$F_{\zeta}^{-1}$ on $[0,\langle\chi,\zeta\rangle)$, which we extend to
be constant and equal to $x^*$ on $[\langle\chi,\zeta\rangle,\infty)$; this
will have a good interpretation for unbounded supports. 

Evidently then,
\begin{equation}\nonumber
  S(t)=F_{\mathcal{B}(t)}^{-1}(t-\beta(t)), \qquad t\ge0,
\end{equation}
which we take as the definition of the cumulative shift of the current batch
$S(t)$. Note that $S(t)$ may become constant equal to $x^*$ for
$t\ge\beta(t)+W(\beta(t))$, as might happen if there are no waiting jobs when
the active batch has completed processing and the server starts an idle
period. Note also that the above considerations still make sense when
$\beta(t)=0$ and $\mathcal{B}_0=0$, as might be the case with zero initial
condition and an initial idle period.  This is the one possibility for which
$\mathcal{B}(t)=0$. In this case $W(\beta(t))=0$,  $x^*=-\infty$,
$F_{\mathcal{B}(t)}$ and $F_{\mathcal{B}(t)}^{-1}$ are identically zero, and
$S(t)=0$.

We can now define 
\begin{equation}\nonumber
  \sigma(t)=\mathcal{B}(t)(\cdot+_0S(t)), \qquad t\ge0,
\end{equation}
which if the profile of residual service times of jobs in the active batch.
Note that if the system is currently idle, $\sigma(t)=0$. The second
part of the state descriptor describing the current service time profile of
waiting jobs is defined
\begin{equation}\nonumber
  \mu(t)=\mathcal{E}(t)-\mathcal{E}(\beta(t)), \qquad t\ge0.
\end{equation}

The pair $(\sigma(\cdot),\mu(\cdot))$ is the state descriptor we use for this
model. It takes values in $\mathbf{D}([0,\infty),\mathcal{M}_1)^2$.  Notice that
the workload process can be written $W(t)=\langle
\chi,\sigma(t)+\mu(t)\rangle$ and the queue length process can be written
$Z(t)=\langle 1,\sigma(t)+\mu(t)\rangle$.

\subsection{Fluid model}\label{s.fluidModel}

We next define a critical fluid model that will serve as the limiting
approximation to the stochastic model. Let $\alpha>0$ and let $\nu$ be a
probability measure on $\mathbb{R}_+$ with $\nu(\{0\})=0$ and finite mean
$\langle\chi,\nu\rangle=1/\alpha$. Let $\xi\in\mathcal{M}_1$ with $\xi(\{0\})=0$
and finite first moment $w=\langle\chi,\xi\rangle<\infty$. Recall that
$\llbracket t\rrbracket_w$ denotes the time since the most recent integer
multiple of $w$.

We define a fluid model path for $\alpha$, $\nu$, and initial condition $\xi$
to be the pair $(\sigma(\cdot),\mu(\cdot))\in
\mathbf{D}([0,\infty),\mathcal{M}_1)^2$  given by
\begin{align}\label{shifting solution}
  \sigma(t) & =
  \begin{cases} \xi_+\bigl( \cdot+F^{-1}_{\xi}
    (t) \bigr), \qquad & t\in[0,w), \\
    \alpha w \nu_+ \left(\cdot +F^{-1}_{\alpha w \nu }
    \left(\llbracket t\rrbracket_w\right)\right), & t\in [w,\infty),
    \end{cases} \\
    \intertext{and}
      \mu(t) & = \alpha \llbracket t\rrbracket_w \nu, \qquad t\ge0.
\label{growing solution}
\end{align}

We refer to $\sigma(\cdot)$ as the shifting path and $\mu(\cdot)$ as the
growing path. Note that if $\xi=0$, both paths are identically zero for all
time. If $\xi\ne0$, then its first moment $w$ is the initial workload in the
fluid model and determines the length of the time intervals over which the
fluid model is periodic. The shifting path has an initial atypical interval
$[0,w)$ during which the initial condition $\xi$ is shifting left according to
the dynamics in the first line of \eqref{shifting solution}. During this time
the growing path is building up in the shape of $\nu$ at rate $\alpha$. At
time $w$, the initial condition has been cleared out and the mass that has
accumulated in the growing path instantaneously becomes the shifting path;
that is $\sigma(w)=\mu(w-)=\alpha w\nu$. Since the fluid model is critical
with $\langle\chi,\nu\rangle=\alpha^{-1}$, the workload at time $w$ also
equals $w$, so the shifting path will empty in another time interval of length
$w$, with the dynamics of $\sigma(\cdot)$ henceforth governed by the second
line in \eqref{shifting solution}. This cycle repeats indefinitely. 

Notice that the above fluid dynamics consist of discontinuous paths (at
integer multiples of $w$), so care is needed in the proof of a limit theorem.
Note also that the fluid model does not converge to a steady state as
$t\to\infty$ but rather, after a finite time $w$, joins a periodic orbit of
states that is determined by the initial workload $w$. This orbiting is not
just the result of separating the state descriptor into the two parts
$\sigma(\cdot)$ and $\mu(\cdot)$. It can readily be seen that for most
distributions $\nu$, the total queue length $\langle
1,\sigma(\cdot)+\mu(\cdot)\rangle$ will also oscillate.


\subsection{Sequence of models and main result}

We now consider a sequence of stochastic models indexed by $r>0$, which tends
to infinity.  Each model is defined as in Section \ref{s.stochasticModel}
under the assumptions stated there, and we append a superscript $r$ to all
symbols associated to the $r$th model. That is, we assume that for each $r$
there are stochastic primitives $\mathcal{E}^r(\cdot)$, and $\mathcal{B}^r_0$
from which are defined performance processes $W^r(\cdot)$, $I^r(\cdot)$, batch
start times and profiles $\{\beta^r_k\}$, $\{\mathcal{B}_k\}$, active batch
processes $\beta^r(\cdot)$, $\ell^r(\cdot)$, $\mathcal{B}^r(\cdot)$, and
$S^r(\cdot)$, and state descriptors $(\sigma^r(\cdot),\mu^r(\cdot))$.  Each
model may be defined on its own probability space
$(\Omega^r,\mathscr{F}^r,P^r)$. 

\textbf{Fluid scaling.} We apply a fluid or law of large numbers scaling to
objects in the sequence as follows. For all $t\ge0$, we define or the primitives,
\begin{equation}
  \nonumber
  \bar{\mathcal{E}}^r(t)=r^{-1}\mathcal{E}^r(rt), \qquad
  \bar{\mathcal{B}}^r_0=r^{-1}\mathcal{B}^r_0,
\end{equation}
for the workload,
\begin{align}
  \nonumber
  \bar{W}^r(t) & =r^{-1}W^r(rt),\\
  \label{e.fluidIdle}
  \bar{I}^r(t) & =
  \sup_{s\in[0,t]}(\bar{W}^r(0)+\langle\chi,\bar{\mathcal{E}}^r(s)\rangle
  -s)^- = r^{-1}I^r(rt), 
\end{align}
and for the state descriptors,
\begin{align*}
  \bar{\sigma}^r(t) &=r^{-1}\sigma^r(rt),\\
  \bar{\mu}^r(t)& =r^{-1}\mu^r(rt).
\end{align*}
Fluid scale versions of the batch start times and profiles can be defined
analogously to \eqref{e.batchStartTimes}--\eqref{e.batchProfiles} by setting
$\bar{\beta}^r_0=0$, and for $k\ge1$ or $t\ge0$,
\begin{align}
  \label{e.fluidBatchStartTimes}
  \bar{\beta}^r_k  & =\inf\{ s \ge
  \bar{\beta}^r_{k-1}+\bar{W}^r(\bar{\beta}^r_{k-1}) 
: \bar{W}^r(s)>0\} 
    = r^{-1}\beta^r_k, \\
  \label{e.fluidBatches}
  \bar{\mathcal{B}}^r_k & = \bar{\mathcal{E}}^r(\bar{\beta}^r_k)
  -\bar{\mathcal{E}}^r(\bar{\beta}^r_{k-1}) 
  = r^{-1}\mathcal{B}^r_k, \\
  \bar{\beta}^r(t) & =\max\{\bar{\beta}^r_k
  : \bar{\beta}^r_k\leq t\}
  = r^{-1}\beta^r(rt), \\
  \bar{\ell}^r(t) & =\max\{j:\bar{\beta}^r_j\le t\}
  = \ell^r(rt), \\
  \bar{\mathcal{B}}^r(t) & = \bar{\mathcal{B}}^r_{\bar{\ell}^r(t)}
  = r^{-1}\mathcal{B}^r(rt),
  \label{e.fluidBatchProfiles}
\end{align}
where either the first or second equality can be taken as the definition, as
they are equivalent.  Note that $\bar{\mathcal{B}}^r_0=\bar{\sigma}^r(0)$.
Then almost surely, for all $t\ge0$, the fluid scaled workload and state
descriptors satisfy
\begin{align}
  \label{e.fluidWDynamics}
  \bar{W}^r(t)
  & =\bar{W}^r(0)+\langle\chi,\bar{\mathcal{E}}^r(t)\rangle -t +\bar{I}^r(t),\\
  \label{e.fluidSigmaDynamics}
  \bar{\sigma}^r(t) 
  & =\bar{\mathcal{B}}^r(t)\left( \cdot+_0
      F^{-1}_{\bar{\mathcal{B}}^r(t)}(t-\bar{\beta}^r(t) \right), \\
      \bar{\mu}^r(t) 
      &  = \bar{\mathcal{E}}^r(t)-\bar{\mathcal{E}}^r(\bar{\beta}^r(t)).
  \label{e.fluidMuDynamics}
\end{align}
In \eqref{e.fluidSigmaDynamics}, we have used the easily verified fact that
$F^{-1}_{r\zeta}(rt)=F^{-1}_\zeta(t)$, for any $\zeta\in\mathcal{M}_1$ and $t\ge0$.

\textbf{Asymptotic assumptions.} We make the following assumptions as
$r\to\infty$. For some probability measure $\nu$ on $\mathbb{R}_+$ with
no atoms and finite mean $\langle\chi,\nu\rangle=1/\alpha$, we have 
\begin{equation}
  \alpha^r\to\alpha, \qquad \nu^r\ws\nu, \qquad
  \langle\chi,\nu^r\rangle\to\langle\chi,\nu\rangle,
  \label{e.critical}
\end{equation}
and moreover that 
\begin{equation}
  \bar{\mathcal{E}}^r(\cdot)\Rightarrow \alpha(\cdot)\nu,
  \label{e.exogenousLimit}
\end{equation}
in $\mathbf{D}([0,\infty),\mathcal{M}_1)$, uniformly on compact time
intervals, where $\alpha(t)=\alpha t$ for all $t\ge0$. Conditions (such as
uniform integrability) under which \eqref{e.exogenousLimit} holds, in
particular guaranteeing that the underlying convergence on $\mathcal{M}_1$ is
in the Wasserstein$_1$ topology, are well known; see for example
\cite{GromPuWi2002}, Lemma A.2. 

For the initial conditions, assume that for some random
$\mathcal{B}_0\in\mathcal{M}_1$ with no atoms almost surely and 
$W_0=\langle\chi,\mathcal{B}_0\rangle<\infty$,
\begin{equation}
  \bar{\mathcal{B}}^r_0\Rightarrow \mathcal{B}_0,
  \qquad \text{and} \qquad
  \mathbb{E}^r[\langle\chi,\bar{\mathcal{B}}^r_0\rangle]
    \to\mathbb{E}[W_0].
  \label{e.initialLimit}
\end{equation}

As is well known, \eqref{e.critical}--\eqref{e.initialLimit}
imply the following classical fluid limit for the workload process.
\begin{prop}
  \label{p.workloadLimit}
 Under the assumptions \eqref{e.critical}--\eqref{e.initialLimit}, the
 sequence of fluid scaled workload processes $\{\bar{W}^r(\cdot)\}$ converges
  in distribution on $\mathbf{D}([0,\infty),\mathbb{R}_+)$ to the process
  $W(\cdot)$ that is almost surely constant and equal to $W_0$.
\end{prop}

The main result of this article is a fluid approximation for the
measure valued state descriptors of the asymptotically critical gated
processor sharing model. 
\begin{thm}\label{t.main}
As $r\to\infty$ under the assumptions
\eqref{e.critical}--\eqref{e.initialLimit}, the sequence of fluid scaled state
descriptors $\{(\bar{\sigma}^r(\cdot), \bar{\mu}^r(\cdot))\}$ converges in
distribution on $\mathbf{D}([0,\infty),\mathcal{M}_1)^2$ to a limit
$(\sigma(\cdot),\mu(\cdot))$ that is almost surely a fluid model path for
$\alpha$, $\nu$, and initial condition $\mathcal{B}_0$.  \end{thm}

The proof is presented in the next three sections, wherein we derive several
preliminary results, establish tightness of the sequence of fluid scaled state
descriptors, and show convergence to the desired limit. The final section then
contains some further discussion application of the result.

\section{Preliminary observations}

In this section we collect several basic facts that will be needed to
prove the limit theorem, including continuity properties of the function
$F_\zeta$, and some basic regularity properties for atom-free measures over
partitions. 

\begin{lem}\label{l.Fcontinuous}
  The mapping $F:\mathcal{M}_1\to\mathbf{C}$ defined by $\zeta\mapsto
F_\zeta(\cdot)$ is  continuous with respect to the topology of uniform
convergence on $\mathbf{C}$.  \end{lem}

\textbf{Proof.} Given $\zeta_n\ws\zeta$ in $\mathcal{M}_1$, we have
$F_{\zeta_n}(x)=\langle\chi_x,\zeta_n\rangle \to \langle \chi_x,\zeta\rangle =
F_\zeta(x)$ for each $x\ge 0$, since $\chi_x\in\mathbf{C}_b$. Since also
$F_{\zeta_n}(\infty)=\langle\chi,\zeta_n\rangle \to \langle\chi,\zeta\rangle =
F_\zeta(\infty)<\infty$, and since $F_{\zeta_n}(\cdot)$ and $F_\zeta(\cdot)$
are continuous and non-decreasing, the convergence is uniform.
\hfill$\blacksquare$

 \begin{lem}\label{l.Finverse}
If  $\zeta_n\ws\zeta\ne0$ and
$x_n\to x<\langle\chi,\zeta\rangle$ in $\mathbb{R}_+$, 
Then
\[  F^{-1}_{\zeta_n}(x_n) \to F^{-1}_\zeta(x).   \]
\end{lem}

\textbf{Proof.} Since $F_{\zeta_n}\to F_\zeta$ uniformly by Lemma
\ref{l.Fcontinuous}, $F_{\zeta_n}$ are eventually continuous and increasing on
$[0,y]$ for some $y\in (x,\langle\chi,\zeta\rangle)$.  So for sufficiently
large $n$, $F^{-1}_{\zeta_n}\to F^{-1}_\zeta$ uniformly on $[0,y]$ which
establishes the result.
\hfill $\blacksquare$\\

For any $\zeta\in\mathcal{M}_1$, let $m_\zeta(x)=\langle
1_{[0,x]},\zeta\rangle$ denote the cumulative mass function of $\zeta$.  


\begin{lem}\label{l.openSets}
If $\zeta\in\mathcal{M}_1$ has no atoms, then for each $\eta>0$ the set
\[
\mathscr{U}_\eta(\zeta)=\{\xi\in\mathcal{M}_1: \|m_\xi-m_\zeta\|_\infty \vee
\|F_\xi-F_\zeta\|_\infty< \eta\} \] 
is open in $\mathcal{M}_1$. 
\end{lem}

\textbf{Proof.} If $\zeta_n\ws\zeta$, then $m_{\zeta_n}\to m_\zeta$ pointwise
since $m_\zeta$ is continuous. So the Glivenko-Cantelli theorem implies that
this convergence is actually uniform. Furthermore, $F_{\zeta_n}\to F_\zeta$
uniformly by Lemma \ref{l.Fcontinuous}. \hfill $\blacksquare$

For $T\in[0,\infty]$, a finite partition of $[0,T)$ is a set of points
$\{x_i\in[0,\infty]:i=0,\ldots,N\}$ such that $0=x_0<x_1<\cdots<x_N$ and
$T\le x_N$. Its mesh
is $\|\{x_i\}\| = \min_{i=1,\ldots N} (x_i-x_{i-1})$, and for
$f:\mathbb{R}_+\to\mathbb{R}$, the modulus of $f$ over the partition is
\[\|f\|_{\{x_i\}}=\max_{i=1,\ldots,N} \sup_{x_{i-1}\le
s<t<x_i}|f(t)-f(s)|.\] Given a partition $\{x_i\}$ of $\mathbb{R}_+$, define
the functions
\begin{equation}
  \psi_i=\chi_{x_i}-\chi_{x_{i-1}}=\chi_{x_i-x_{i-1}}(\cdot-x_{i-1}),
  \qquad i=1,\ldots,N,
\label{e.chiIncrements}
\end{equation}
which are the increments of the identity function $\chi$ over the partition. 

\begin{lem}\label{l.goodPartitions}
  Let $\zeta\in\mathcal{M}_1$ have no atoms. Then for all $\varepsilon>0$
  there exists a finite partition $\{x_i\}$ of $\mathbb{R}_+$, such that
  \[ \|m_\zeta\|_{\{x_i\}} \vee \|F_\zeta\|_{\{x_i\}}\le\varepsilon,\]
  and $\min\{\langle\psi_i,\zeta\rangle :i=1,\ldots,N\}>0$.
\end{lem}

\textbf{Proof.} Let $x_0=0$ and for $i\ge 1$ inductively define
$x_i=\sup\{x>x_{i-1}: (m_\zeta(x)-m_\zeta(x_{i-1})) \vee 
  (F_\zeta(x)-F_\zeta(x_{i-1}))\le\varepsilon$. Then $x_N=\infty$ for some
  finite $N$ since $\zeta$ has finite total mass and first moment. The
  intervals $[x_{i-1},x_i)$ are nonempty for $i=1,\ldots,N$, by continuity of
  $m_\zeta$ and $F_\zeta$, and so $\{x_i\}_{i=1}^N$ is a finite partition of
  $\mathbb{R}_+$. Lastly, continuity also implies that it is impossible for
  $x_i$ to equal the supremum of the support of $\zeta$ for $i<N$ (if it did,
  it would already have to equal infinity). Thus the supports of all
  $\psi_i$, $i=1,\ldots,N$, intersect the support of $\zeta$ non-trivially,
  which yields the final property. \hfill $\blacksquare$

\section{Tightness of the state descriptors}

We now show that the sequences of fluid scaled state
descriptors $\{\bar{\sigma}^r(\cdot)\}$ and $\{\bar{\mu}^r(\cdot)\}$ are tight. 

\textbf{Tightness of $\{\bar{\sigma}^r(\cdot)\}$}. 
  Let $\mathbf{C}_c\subset\mathbf{C}_b$ be
the functions with compact support.  Then $\mathbf{C}_c$ separates points of
$\mathcal{M}_1$, is closed under addition, and $\zeta\mapsto \langle
f,\zeta\rangle$ is continuous on $\mathcal{M}_1$ for each $f\in\mathbf{C}_c$. So
by Jakubowski's tightness criterion (see \cite{Jaku1986} Theorem 3.1), it
suffices to show that for each $T>0$, $\bar{\sigma}^r(\cdot)$ is compactly
contained on $[0,T]$ with high probability, and $\{\langle
f,\bar{\sigma}^r(\cdot)\rangle\}$ is tight in $\mathbf{D}([0,T],\mathbb{R}_+)$
for each $f\in\mathbf{C}_c$. 

\textbf{Compact containment.} Let $\varepsilon>0$ and observe that by
\eqref{e.exogenousLimit}, $\{\bar{\mathcal{E}}^r(T)\}$ is tight: there exists
a compact $\mathcal{K}\subset\mathcal{M}_1$, containing $\alpha T\nu$, such that
\[
\liminf_{r\to\infty}P^r(\bar{\mathcal{E}}^r(T)\in\mathcal{K})\ge
1-\varepsilon.
\]
Let $M=\sup_{\zeta\in\mathcal{K}}\langle 1,\zeta\rangle<\infty$ and define
$G_x(\zeta)= \langle\chi(\cdot-x),\zeta\rangle = \langle\chi,\zeta\rangle -
\langle \chi_x,\zeta\rangle$ for $x\in\mathbb{R}_+$.  Then $G_x$ is continuous
on $\mathcal{M}_1$ so
$\varphi(x)=\sup_{\zeta\in\mathcal{K}}G_x(\zeta)<\infty$. Moreover, monotone
convergence implies $G_x(\zeta)\to 0$ as $x\to \infty$. That is, $G_x(\cdot)$
decreases pointwise monitonically to zero on $\mathcal{K}$, and thus
uniformly. So $\varphi(x)\to 0$ as $x\to\infty$.  The set 
\[ \tilde{\mathcal{K}}=\{\zeta\in\mathcal{M}_1:\langle 1,\zeta\rangle\le
M \text{ and } \langle\chi(\cdot-x),\zeta\rangle \le \varphi(x)
  \text{ for all $x$}\}\]
  is easily seen to be compact in $\mathcal{M}_1$. So since 
  \eqref{e.fluidSigmaDynamics},
\eqref{e.fluidBatchProfiles}, and \eqref{e.fluidBatches} imply the bounds
\begin{align}\label{e.massBound}
  \langle 1,\bar\sigma^r(t)\rangle 
   &\le \langle 1,\bar{\mathcal{B}}^r(t)\rangle
   \le \langle 1,\bar{\mathcal{E}}^r(t)\rangle, & \quad t\in [0,T], \\
   \langle\chi(\cdot -x),\bar{\sigma}^r(t)\rangle
   &\le\langle\chi(\cdot-x),\bar{\mathcal{B}}^r(t)\rangle
   \le \langle\chi(\cdot -x),\bar{\mathcal{E}}^r(T)\rangle, & \quad t\in[0,T],
  \label{e.workloadBound}
\end{align}
we have $\liminf_{r\to\infty}P^r(\bar{\sigma}^r(t)\in\tilde{\mathcal{K}}
\text{ for all $t\in[0,T]$})\ge1-\varepsilon$, which establishes compact
containment.

\textbf{Tightness of projections.} 
Next, we fix $T>0$ and a nonzero $f\in\mathbf{C}_c$. Showing tightness of the
real-valued processes  $\{\langle f,\bar{\sigma}^r(\cdot)\rangle\}$ in
$\mathbf{D}([0,T],\mathbb{R}_+)$ requires a compact containment condition, and
uniform control with high probability of the $J_1$-modulus of continuity 
\[
\operatorname{w}'(z(\cdot),\delta,T)
= \inf_{\{t_i\}} \max_{i=1,\ldots,N}\sup_{s,t\in[t_{i-1},t_i\wedge T)} |z(t)-z(s)|,
\]
where the infimum is over all finite partitions $\{t_i\}$ of $[0,T)$ with mesh
$\|\{t_i\}\|\ge \delta$; see [?]. Since $\langle f,\bar{\sigma}^r(t)\rangle\le
\|f\|_\infty\langle 1,\bar{\sigma}^r(t)\rangle$ for all $t\ge0$, the
real-valued compact containment follows immediately from the measure-valued
compact containment already established above.
So it remains to show that for all $\varepsilon>0$ there exists
$\delta_\varepsilon>0$ such that 
\begin{equation}\label{e.modulus}
  \liminf_{r\to\infty}P^r\left( \operatorname{w}'\left(\langle
  f,\bar{\sigma}^r(\cdot)\rangle ,\delta_\varepsilon,T\right)
 \le\varepsilon \right)\ge 1-\varepsilon.
\end{equation}

To that end, let $\varepsilon>0$ and define a series of constants from
$\varepsilon$, $f$, $T$, $\alpha$, and $\nu$ as follows. First, since $f$ is
uniformly continuous, we can choose $0<\delta_f<1$ such that \[
\sup_{a\le\delta_f} \|f-f(\cdot-a)\|_\infty \le\varepsilon/4\alpha T. \]
Next choose $0<\varepsilon_0<\varepsilon\delta_f/8\|f\|_\infty$. Since $\nu$ has no atom at
zero, choose $\theta>0$  such that $\langle 1_{[0,\theta)},\alpha T\nu\rangle
\le \varepsilon_0$, set $c=\varepsilon_0\theta$ and choose $0<b<c$.

Next, since $\alpha T\nu$ has no atoms, use Lemma \ref{l.goodPartitions} to 
choose a finite partition $\{x_i\}_{i=0}^N$ of $\mathbb{R}_+$ such that
\[ \|m_{\alpha T\nu}\|_{\{x_i\}} \vee \|F_{\alpha T\nu}\|_{\{x_i\}}
\le\varepsilon_0,\] 
and $ \min\{\langle\psi_i,\alpha T\nu\rangle:i=1,\ldots,N\}>0$,
where $\psi_i$ are the increments of $\chi$ over the partition, as in
\eqref{e.chiIncrements}. It follows that also
\[ \delta_1 = \min\{\langle\psi_i,\alpha b\nu\rangle:i=1,\ldots,N\}>0.\]

We define a set of partitions for the initial condition in a similar way. By
\eqref{e.initialLimit} there is a compact $\mathcal{L}\subset\mathcal{M}_1$
such that 
\begin{equation}
  \liminf_{r\to\infty}P^r\left( \bar{\mathcal{B}}^r_0\in\mathcal{L}\right)\ge
  1-\varepsilon.
  \label{e.initCompact}
\end{equation}
Note that if $\mathcal{A}$ is a full probability set containing atom-free
measures, then the closure of $\mathcal{K}\cap\mathcal{A}$ is still compact
with probability equal to that of $\mathcal{K}$. So since $\mathcal{B}_0$
almost surely has no atoms, we can assume without loss of generality that the
elements of $\mathcal{K}$ are atom-free.  By Lemma \ref{l.openSets}, we can
choose a finite set $\{\xi^j\in\mathcal{L}:j=1,\ldots,L\}$ such that the open
neighborhoods $\{\mathscr{U}_\eta(\xi^j):j=1,\ldots,L\}$ cover $\mathcal{L}$.
For each $j=1,\ldots,L$, apply Lemma \ref{l.goodPartitions} again to choose a
finite partition $\{y^j_i\}_{i=0}^{N_j}$ of $\mathbb{R}_+$ such that
$\|m_{\xi^j}\|_{\{y^j_i\}} \vee \|F_{\xi^j}\|_{\{y^j_i\}} \le \varepsilon_0$
and  
  \[\delta_0=\min\{\langle\psi^j_i,\xi^j\rangle:j=1,\ldots,L,\quad i=1,\ldots,N_j\} 
  >0 .\]
where $\psi^j_i=\chi_{y^j_i-y^j_{i-1}}(\cdot-y^j_{i-1})$ are the increments of
$\chi$ over $\{y^j_i\}$ for each $j=1,\ldots,L$.

Finally let $\delta_\varepsilon=(\delta_0\wedge\delta_1)/2$ and 
choose $\eta<(\varepsilon_0/2) \wedge (\delta_\varepsilon/4)$. 
Having set up the various constants, we now establish the 
events on which \eqref{e.modulus} holds.  The real-valued functions
$\zeta\mapsto \langle 1_{[0,x_i]},\zeta\rangle = m_\zeta(x_i)$ are continuous
at each point $\alpha t\nu\in\mathcal{M}_1$, $t\ge0$, because $\nu$ has no atoms.
Therefore \eqref{e.exogenousLimit} implies that
$m_{\bar{\mathcal{E}}^r(\cdot)}(x_i)$ converges uniformly on $[0,T]$ to
$m_{\alpha(\cdot)\nu}(x_i)$ for each $i=0,\ldots,N$. Similarly,
$F_{\bar{\mathcal{E}}^r(\cdot)}(x_i)$ converges uniformly on $[0,T]$ to
$F_{\alpha(\cdot)\nu}(x_i)$ for each $i$, since $\chi_{x_i}\in\mathbf{C}_b$
for $i<N$ and $\chi_\infty=\chi$.
So noting that $m_\zeta-m_\xi=m_{\zeta-\xi}$ whenever
$\zeta-\xi\in\mathcal{M}_1$, and letting 
\begin{multline}
  \nonumber
  \Omega^r_1=\left\{\max_{i=0,\ldots,N}\sup_{0\le s<t<T}\left|
m_{\bar{\mathcal{E}}^r(t)-\bar{\mathcal{E}}^r(s)}(x_i) -
m_{\alpha(t-s)\nu}(x_i)\right| \right. \\
\vee \left. \left| F_{\bar{\mathcal{E}}^r(t)-\bar{\mathcal{E}}^r(s)}(x_i) -
F_{\alpha(t-s)\nu}(x_i)\right| \le \eta   \right\},
\end{multline}
we see that $\liminf_{r\to\infty}P^r\left( \Omega^r_1 \right)=1.$ Letting
\begin{equation}
    \Omega^r_b=\left\{\inf_{t\in[0,T]}\bar{W}^r(t)\ge b\right\},\qquad
    \Omega^r_c=\left\{\sup_{t\in[0,T]}\bar{W}^r(t)\le c\right\},
  \nonumber
\end{equation}
Proposition \ref{p.workloadLimit} implies that $\liminf_{r\to\infty}P^r\left(
\Omega^r_b\cup\Omega^r_c \right)=1.$ Thus letting $\Omega^r_0$ denote the
events in \eqref{e.initCompact} and $\Omega^r_*=\Omega^r_0\cap\Omega^r_1\cap
(\Omega^r_b\cup\Omega^r_c)$, we have $\liminf_{r\to\infty}\left( \Omega^r_*
\right)=1$, and it suffices to show that on $\Omega^r_*$,
\begin{equation}
  \operatorname{w}'(\langle f,\bar{\sigma}^r(\cdot)\rangle,\delta_\varepsilon,T)
    \le \varepsilon.
    \label{e.modulusToShow}
\end{equation}

For this, fix $\omega\in\Omega^r_*$ and suppose first that
$\omega\in\Omega^r_b$. We define a partition of $[0,T)$ for this particular
sample path as follows. Let $\xi=\bar{\sigma}^r(0)$ and since
$\xi\in\mathcal{L}$, choose a neighborhood $\mathscr{U}_\eta(\xi^j)$
containing $\xi$.  Using the partition of $\mathbb{R}_+$ corresponding to
$\xi^j$, let $t^0_i=F_{\xi}(y^j_i)$ for $i=0,\ldots,N_j-1$ and
$t^0_{N_j}=\bar{\beta}^r_1$. For $k=1,\ldots,\bar{\ell}^r(T)$, let
$t^k_i=\bar{\beta}^r_k+F_{\bar{\mathcal{B}}^r_k}(x_i)$ for $i=0,\ldots,N-1$
and $t^k_N=\bar{\beta}^r_{k+1}$.  Note that $t^0_{N_j}=\bar{\beta}^r_1=t^1_0$
and $t^k_N=\bar{\beta}^r_{k+1}=t^{k+1}_0$ for $k=1,\ldots\bar{\ell}^r(t)-1$,
and we consider these unique points with multiple labels. Then $\{t_i\}=
\bigcup_{k=0,\ldots,\bar{\ell}^r(T)} \{t^k_i\}$ is a finite partition of
$[0,T)$.  To compute its mesh, observe that for $k=1,\ldots,\bar{\ell}^r(T)$
and $i=1,\ldots,N$,
\begin{multline}
  t^k_{i}-t^k_{i-1} =
  F_{\bar{\mathcal{B}}^r_k}(x_i)-F_{\bar{\mathcal{B}}^r_k}(x_{i-1})\\
  = F_{\bar{\mathcal{E}}^r(\bar{\beta}^r_k)-\bar{\mathcal{E}}^r(\bar{\beta}^r_{k-1})}(x_i)
      -F_{\bar{\mathcal{E}}^r(\bar{\beta}^r_k)
           -\bar{\mathcal{E}}^r(\bar{\beta}^r_{k-1})}(x_{i-1}) \\
  \ge F_{\alpha(\bar{\beta}^r_k-\bar{\beta}^r_{k-1})\nu}(x_i)
      -F_{\alpha(\bar{\beta}^r_k-\bar{\beta}^r_{k-1})\nu}(x_{i-1}) -2\eta\\
      = \left\langle\chi_{x_i-x_{i-1}}(\cdot-x_{i-1}),
      \alpha(\bar{\beta}^r_k-\bar{\beta}^r_{k-1})\nu\right\rangle - 2\eta \\
      \ge \left\langle\psi_i, \alpha b\nu\right\rangle - 2\eta 
      \ge \delta_1-2\eta \ge \delta_\varepsilon,
  \nonumber
\end{multline}
since $\bar{\beta}^r_k\ge\bar{\beta}^r_{k-1}+b$ on $\Omega^r_b$. For
$k=0$, the estimate is derived in the same way, with $\{y^j_i\}$ in place of
$\{x_i\}$, $\bar{\mathcal{B}}^r_0=\xi$ in place of $\bar{\mathcal{B}}^r_k$,
$\xi^j$ in place of $\alpha(\bar{\beta}^r_k-\bar{\beta}^r_{k-1})\nu$, and
$\delta_0$ in place of $\delta_1$. It follows that 
$\|\{t_i\}\|\ge\delta_\varepsilon$. 

Let $t^k_{i-1}\le s<t<t^k_i\wedge T$ for some $k=1,\ldots,\bar{\ell}^r(t)$ and
$i=1,\ldots,N$. If it happens that $i=N$ and that
$\bar{\beta}^r_k+F_{\bar{\mathcal{B}}^r_k}(x_N)<t^k_N=\bar{\beta}^r_{k+1}$,
then it is because the system emptied at $t^k_N$ and is idle until the next
start time. In this case $\bar{\sigma}^r(\cdot)$ is identically zero during
$[t^k_N,\bar{\beta}^r_{k+1})$, so we may assume without loss of generality in
this case that $t<\bar{\beta}^r_k+F_{\bar{\mathcal{B}}^r_k}(x_N)$. Then
$F_{\bar{\mathcal{B}}^r_k}(x_{i-1})\le
s-\bar{\beta}^r_k<t-\bar{\beta}^r_k<F_{\bar{\mathcal{B}}^r_k}(x_i)$, and
so $x_{i-1}\le F^{-1}_{\bar{\mathcal{B}}^r_k}(s-\bar{\beta}^r_k) <
F^{-1}_{\bar{\mathcal{B}}^r_k}(t-\bar{\beta}^r_k) < x_i$. This implies
\begin{align*}
    m_{\bar{\mathcal{B}}^r_k}(\bar{S}^r(t))
  -m_{\bar{\mathcal{B}}^r_k}(\bar{S}^r(s))
   &\le m_{\bar{\mathcal{B}}^r_k}(x_i)
   -m_{\bar{\mathcal{B}}^r_k}(x_{i-1}) \\
 & \le  m_{\alpha(\bar{\beta}^r_k- \bar{\beta}^r_{k-1})\nu}(x_i)
-m_{\alpha(\bar{\beta}^r_k-\bar{\beta}^r_{k-1})\nu}(x_{i-1})+2\eta \\
  &\le  m_{\alpha T\nu}(x_i)
-m_{\alpha T\nu}(x_{i-1})+2\eta \\
  & \le \varepsilon_0 +2\eta\le2\varepsilon_0,
 \end{align*}
 and also
\begin{align*}
    F_{\bar{\mathcal{B}}^r_k}(\bar{S}^r(t))
  -F_{\bar{\mathcal{B}}^r_k}(\bar{S}^r(s))
   &\le F_{\bar{\mathcal{B}}^r_k}(x_i)
   -F_{\bar{\mathcal{B}}^r_k}(x_{i-1}) \\
 & \le  F_{\alpha(\bar{\beta}^r_k- \bar{\beta}^r_{k-1})\nu}(x_i)
-F_{\alpha(\bar{\beta}^r_k-\bar{\beta}^r_{k-1})\nu}(x_{i-1})+2\eta \\
  &\le  F_{\alpha T\nu}(x_i)
-F_{\alpha T\nu}(x_{i-1})+2\eta \\
  & \le \varepsilon_0 +2\eta\le2\varepsilon_0,
 \end{align*}
by definition of $\Omega^r_1$ and $\{x_i\}$.  Applying this to the bound
\begin{equation}\label{e.oscBound}
\begin{aligned}
|\langle f,\bar{\sigma}^r(t)\rangle -&\langle f,\bar{\sigma}^r(s)\rangle| 
  = \left|\left\langle f,\bar{\mathcal{B}}^r_k\left(\cdot +_0
\bar{S}^r(t)\right)\right\rangle 
 - \left\langle f,\bar{\mathcal{B}}^r_k\left(\cdot +_0
 \bar{S}^r(s)\right)\right\rangle\right|\\
 & \le\left\langle\left| (1_{(0,\infty)}f)(\cdot 
-\bar{S}^r(t))
        -(1_{(0,\infty)}f)(\cdot- \bar{S}^r(s))\right|
        ,\bar{\mathcal{B}}^r_k \right\rangle\\
     &   \le
        \|f\|_\infty\langle 1_{(\bar{S}^r(s), \bar{S}^r(t)]},
         \bar{\mathcal{B}}^r_k\rangle \\
     &\quad    + \left\|f-f\left(\cdot-(\bar{S}^r(t)-\bar{S}^r(s))\right)\right\|_\infty
     \langle 1_{(\bar{S}^r(t),\infty)},\bar{\mathcal{B}}^r_k\rangle,
\end{aligned}
\end{equation}
we see that the first right hand term is bounded above by
$\|f\|_\infty2\varepsilon_0\le\varepsilon/2$. If $\bar{S}^r(t)-\bar{S}^r(s)\le\delta_f$, bound
the second right hand term by
\[   \sup_{a\le \delta_f}\|f - f(\cdot-a))\|_\infty \bar{E}^r(T)
\le \frac{\varepsilon}{4\alpha T}2\alpha T = \frac{\varepsilon}{2}.
  \]
If instead $\bar{S}^r(t)-\bar{S}^r(s)>\delta_f$, bound this term by 
\begin{align*}
 2\|f\|_\infty
 \frac{\langle\chi_{\bar{S}^r(t)}-\chi_{\bar{S}^r(s)},\bar{\mathcal{B}}^r_k\rangle}
 {\bar{S}^r(t)-\bar{S}^r(s)}
 & \le
 \frac{2\|f\|_\infty}{\delta_f}(F_{\bar{\mathcal{B}}^r_k}(\bar{S}^r(t)) 
 -F_{\bar{\mathcal{B}}^r_k}(\bar{S}^r(s))) \\
 & =\frac{2\|f\|_\infty}{\delta_f}2\varepsilon_0
    \le\frac{\varepsilon}{2}.
 \end{align*}
 In either case, we see that \eqref{e.oscBound} implies that $|\langle
 f,\bar{\sigma}^r(t)\rangle -\langle f,\bar{\sigma}^r(s)\rangle|\le
 \varepsilon$.  For the case $t^0_{i-1}\le s<t<t^0_i \wedge T$, the above
 argument is identical on $[0,\bar{\beta}^r_1)$ with
 $\bar{\mathcal{B}}^r_0=\xi$ in place of $\bar{\mathcal{B}}^r_k$, $\{y^j_i\}$
 in place of $\{x_i\}$, and $\xi^j$ in place of $\alpha T\nu$. Therefore, we
 have established \eqref{e.modulusToShow}.

\section{Convergence to the limit}

Having established that $\{\bar{\sigma}^r(\cdot)\}$ and
$\{\bar{\mu}^r(\cdot)\}$ are tight, and using the convergence in
\eqref{e.exogenousLimit}, \eqref{e.initialLimit}, and Propostion
\ref{p.workloadLimit}, we conclude that 
\[ \{(\bar{\mathcal{E}}^r(\cdot), \bar{W}^r(\cdot),\bar{\sigma}^r(\cdot),
\bar{\mu}^r(\cdot))\} \] 
is jointly tight. By passing to a subsequence, we get convergence in
distribution 
\begin{equation}
  (\bar{\mathcal{E}}^r(\cdot), \bar{W}^r(\cdot),\bar{\sigma}^r(\cdot),
\bar{\mu}^r(\cdot)) 
\Rightarrow
(\alpha(\cdot)\nu,W^*(\cdot),\sigma^*(\cdot),\mu^*(\cdot)),
\label{e.dconv}
\end{equation}
where almost surely
$W^*(\cdot)=\langle\chi,\sigma^*(\cdot)+\mu^*(\cdot)\rangle$ is constant and
$\sigma^*(0)=\xi$, a random measure with $\xi(\{0\})=0$ and
$\langle\chi,\xi\rangle=W^*(0)$. 

To prove Theorem \ref{t.main}, we must uniquely characterize the last two
components in the limit above, for which it suffices to prove that almost
surely, they are a fluid model path for $\alpha$, $\nu$, and initial condition
$\xi$ as defined in Section \ref{s.fluidModel}. 

By the Skorohod representation theorem, there exists a sequence of tuples
$(\tilde{\mathcal{E}}^r(\cdot), \tilde{W}^r(\cdot),\tilde{\sigma}^r(\cdot),
\tilde{\mu}^r(\cdot))$ that are defind on a common probability space
$(\Omega,\mathscr{F},P)$ and are equal in distribution to the tuples on the
left side of \eqref{e.dconv}, and there exists a tuple
$(\alpha(\cdot)\nu,W(\cdot),\sigma(\cdot),\mu(\cdot))$ on $\Omega$ that is
equal in distribution to the right side of \eqref{e.dconv} such that almost
surely  on $\Omega$,
\begin{equation}
  (\tilde{\mathcal{E}}^r(\cdot),
\tilde{W}^r(\cdot),\tilde{\sigma}^r(\cdot),
\tilde{\mu}^r(\cdot))
\rightarrow
(\alpha(\cdot)\nu,W(\cdot),\sigma(\cdot),\mu(\cdot)),
\label{e.skorohodConv}
\end{equation}
uniformly on compact time intervals in the first two components, and in the
Skorohod $J_1$-topology in the second two. 

For each $r$, we construct batch start times and profiles for the Skorohod
representations in the same way as
\eqref{e.fluidBatchStartTimes}--\eqref{e.fluidBatchProfiles} for the original
fluid scaled models.  That is, $\tilde{\beta}^r_0=0$ and
$\tilde{\mathcal{B}}^r_0=\tilde{\sigma}^r(0)$, and for $k\ge1$ or $t\ge0$,
\begin{align*}
  \tilde{\beta}^r_k  & =\inf\{ s \ge
  \tilde{\beta}^r_{k-1}+\tilde{W}^r(\tilde{\beta}^r_{k-1}) 
  : \tilde{W}^r(s)>0\}, \\
  \tilde{\mathcal{B}}^r_k & = \tilde{\mathcal{E}}^r(\tilde{\beta}^r_k)
  -\tilde{\mathcal{E}}^r(\tilde{\beta}^r_{k-1}), \\
  \tilde{\beta}^r(t) & =\max\{\tilde{\beta}^r_k
  : \tilde{\beta}^r_k\leq t\}, \\
  \tilde{\ell}^r(t) & =\max\{j:\tilde{\beta}^r_j\le t\}, \\
  \tilde{\mathcal{B}}^r(t) & = \tilde{\mathcal{B}}^r_{\tilde{\ell}^r(t)}.
\end{align*}
As before we define
$\tilde{S}^r(t)=F^{-1}_{\tilde{\mathcal{B}}^r(t)}(t-\tilde{\beta}^r(t))$ and we can
also define $\tilde{I}^r(t)=\sup_{s\in[0,t]}(\tilde{W}^r(0)+\langle\chi ,
\tilde{\mathcal{E}}^r(s)\rangle -s)^-$ analogously to \eqref{e.fluidIdle}.

It is not difficult, though somewhat tedious, to verify that these are
measurable functions of the tuples in \eqref{e.skorohodConv}. Since they are
the same functions used in \eqref{e.fluidIdle}--\eqref{e.fluidBatchProfiles},
we see that the above objects have the same distributions as their original
counterparts and that consequently, $\tilde{W}^r(\cdot)$,
$\tilde{\sigma}^r(\cdot)$, and $\tilde{\mu}^r(\cdot)$ satisfy the same
properties almost surely as do their original analogs. That is, almost surely
for all $t\ge0$,
\begin{equation}
  \label{e.skorohodDynamics}
\begin{aligned}
  \tilde{W}^r(t)&=\tilde{W}^r(0)+\langle\chi,\tilde{\mathcal{E}}^r(t)\rangle
  -t +\tilde{I}^r(t), \\
  \tilde{\sigma}^r(t)&=\tilde{\mathcal{B}}^r(t)
  \left( \cdot+_0F^{-1}_{\tilde{\mathcal{B}}^r(t)}(t-\tilde{\beta}^r(t))
  \right),\\
  \tilde{\mu}^r(t)&=\tilde{\mathcal{E}}^r(t)-\tilde{\mathcal{E}}^r(\tilde{\beta}^r(t)).
\end{aligned}
\end{equation}

\subsection{Zero initial condition}\label{s.zero}

We first consider the case $W(0)=0$. This will use the following simple bound on
the total mass of the shifting path. 
\begin{lem}
\label{l.zeroBound}
Almost surely, for all $r$ and $t\in[0,T]$,
  \begin{equation}
    \nonumber
    \langle 1,\tilde{\sigma}^r(t)\rangle
    \le \begin{cases} \langle 1,\tilde{\sigma}^r(0)\rangle,
      & t< \tilde{\beta}^r_1, \\
      \displaystyle\max_{1\le k\le\tilde{\ell}^r(t)}
      \tilde{E}^r\left(\tilde{\beta}^r_k\right)
    -\tilde{E}^r\left(\tilde{\beta}^r_{k-1}\right),
    & t\ge \tilde{\beta}^r_1,
  \end{cases}
\end{equation}
where $\tilde{E}^r(t)=\langle 1,\tilde{\mathcal{E}}^r(t)\rangle$ is the total
mass.
\end{lem}

\textbf{Proof.} If $t<\tilde{\beta}^r_1$, then $\tilde{\beta}^r(t)=0$ and so
$\tilde{\mathcal{B}}^r(t)=\tilde{\mathcal{B}}^r_0=\tilde{\sigma}^r(0)$. The
bound then follows from \eqref{e.skorohodDynamics}. If
$t\ge\tilde{\beta}^r_1$, then $\tilde{\ell}^r(t)\ge1$ and
$\tilde{\mathcal{B}}^r(t)=\tilde{\mathcal{E}}^r(\tilde{\beta}^r_{\tilde{\ell}^r(t)})
- \tilde{\mathcal{E}}^r(\tilde{\beta}^r_{\tilde{\ell}^r(t)-1})$ and the bound
follows from \eqref{e.skorohodDynamics}. \hfill$\blacksquare$

We next show that when the limiting initial condition is zero, intervals
between prelimit batch start times become uniformly small on compact time
intervals.
\begin{lem}\label{l.zeroBetasBunch}
  Almost surely on the event $W(0)=0$, 
  \[
    \lim_{r\to\infty} 
    \max_{1\le k\le\tilde{\ell}^r(T)}(\tilde{\beta}^r_k-\tilde{\beta}^r_{k-1})
    =0,
  \]
for all $T>0.$
\end{lem}
\textbf{Proof.} Fix $\omega\in\Omega$ such that $W(0)=0$ and
\eqref{e.skorohodConv} holds, and let $T>0$. Then uniformly on $[0,T]$,
$\tilde{W}^r(\cdot)\to0$ and also
$\langle\chi_1,\tilde{\mathcal{E}}^r(\cdot)\rangle \to \alpha(\cdot)\langle
\chi_1,\nu\rangle$,  since $\chi_1=\chi\wedge 1$ is bounded and continuous.
Given $\varepsilon>0$, choose
$\delta<\alpha\varepsilon\langle\chi_1,\nu\rangle/2$. Then for sufficiently
large $r$, we have
\[
  \langle\chi_1,\tilde{\mathcal{E}}^r(t)\rangle
  -\langle\chi_1,\tilde{\mathcal{E}}^r(s)\rangle
  \ge \alpha(t-s)\langle\chi_1,\nu\rangle -2\delta,
\]
for all $s\le t\le T$, and we have
\[
  \sup_{t\in[0,T]}\tilde{W}^r(t)<\varepsilon.
\]
Consider any $1\le k\le\tilde{\ell}^r(T)$. Let
$u=\tilde{\beta}^r_{k-1}+\varepsilon$ and suppose that $\tilde{\beta}^r_k >
u$. Then $\tilde{\beta}^r(u)=\tilde{\beta}^r_{k-1}$ and so 
\begin{multline*}
  \tilde{W}^r(u)
     = \langle\chi,\tilde{\sigma}^r(u)+\tilde{\mu}^r(u)\rangle 
     \ge \langle\chi_1,\tilde{\mu}^r(u)\rangle \\
     = \langle\chi_1,\tilde{\mathcal{E}}^r(u)
     -\tilde{\mathcal{E}}^r(\tilde{\beta}^r_{k-1})\rangle 
     \ge \alpha(u-\tilde{\beta}^r_{k-1})\langle\chi_1,\nu\rangle-2\delta
     >0.
\end{multline*}
Since $u=\tilde{\beta}^r_{k-1} + \varepsilon
>\tilde{\beta}^r_{k-1}+ \tilde{W}^r(\tilde{\beta}^r_{k-1})$, this implies $u$
is such a time with positive workload, which contradicts the definition of
$\tilde{\beta}^r_k$. We conclude that $\tilde{\beta}^r_k\le u$ and thus
$\tilde{\beta}^r_k -\tilde{\beta}^r_{k-1}\le\varepsilon$.
\hfill$\blacksquare$

We now show that almost surely, $\sigma(\cdot)$ and $\mu(\cdot)$ are
identically zero on the event $W(0)=0$. Fix $\omega$ on this event such that 
\eqref{e.skorohodConv} and Lemma \ref{l.zeroBetasBunch} holds. 
Let $T>0$ and let $t\in[0,T]$ be any
continuity point of $\sigma(\cdot)$. Then 
$\tilde{\sigma}^r(t)\ws\sigma(t)$. So by Lemma
\ref{l.zeroBound},
\begin{align*}
  \langle 1,\sigma(t)\rangle & = \lim_{r\to\infty} \langle
  1,\tilde{\sigma}^r(t)\rangle \\
 & \le \lim_{r\to\infty}\left(\langle 1,\tilde{\sigma}^r(0)\rangle
  \vee \max_{1\le k\le\tilde{\ell}^r(t)}
  \left(\tilde{E}^r(\tilde{\beta}^r_k)-
\tilde{E}^r(\tilde{\beta}^r_{k-1})\right)\right),
\end{align*}
where we take the maximum to be zero when $\tilde{\ell}^r(t)=0$. By
\eqref{e.skorohodConv}, both
$\tilde{\sigma}^r(0)\ws\sigma(0)$ and
$\tilde{E}^r(\cdot)\to\alpha(\cdot)$ uniformly on $[0,T]$. The first
convergence implies that the first term above vanishes, and the second
convergence combined with Lemma \ref{l.zeroBetasBunch} implies that the second
term vanishes as well. Thus $\sigma(t)=0$ for all continuity points, and
therefore for all $t\in[0,T]$ by right-continuity. 

Similarly, let $t\in[0,T]$ by any continuity point of $\mu(\cdot)$. Then
$\tilde{\mu}^r(t)\ws\mu(t)$ and
\begin{align*}
  \langle 1,\mu(t)\rangle & = \lim_{r\to\infty}
         \langle 1,\tilde{\mu}^r(t)\rangle \\
  & = \lim_{r\to\infty}\left(\tilde{E}^r(t) 
       -\tilde{E}^r(\tilde{\beta}^r(t))\right) \\
      & \le \lim_{r\to\infty}\left(\max_{1\le k\le\tilde{\ell}^r(T)}
         \tilde{E}^r(\tilde{\beta}^r_k) 
       -\tilde{E}^r(\tilde{\beta}^r_{k-1})\right).
\end{align*}
This limit is zero by Lemma \ref{l.zeroBetasBunch} and the uniform convergence
of $\tilde{E}^r(\cdot)$. Again by right-continuity, $\mu(\cdot)$ is
identically zero on $[0,T]$.

\subsection{Nonzero initial condition}\label{s.nonzero}

We next consider the case $\langle\chi,\xi\rangle =W(0)>0$, and show that on this event
$(\sigma(\cdot),\mu(\cdot))$ is a fluid model path for $\alpha$, $\nu$ and
initial condition $\xi$. We first show that in this case almost
surely, batch start times converge to integer multiples of $W(0)$. 

\begin{lem}\label{l.betas}
Almost surely on the event $W(0)>0$, $\tilde{\beta}^r_k \to kw$ for all $T>0$
and all integers, $0\le k \le \left\lfloor T/w\right\rfloor.$  \end{lem}

\textbf{Proof.}  Fix $\omega \in \Omega$ such that $w=W(0)>0$ and
\eqref{e.skorohodConv} holds. Let $T>0$.  Clearly $\tilde{\beta}^r_0\to0.$
Proceeding by induction, assume that $\tilde{\beta}^r_k\to kw$ for some $0\le
k<\lfloor T/w\rfloor$. By \eqref{e.skorohodConv} and since $W(\cdot)$ is
constant and equal to $w>0$, we have that for sufficiently large $r$,
$\inf_{t\in[0,T+1]} \tilde{W}^r(t)>0.$ This implies that for large $r$,
$\tilde{\beta}^r_k+\tilde{W}^r(\tilde{\beta}^r_k)<T+1$. Then since the
workload is not zero in $[0,T+1]$ we have by definition
$\tilde{\beta}^r_{k+1}=\tilde{\beta}^r_k+\tilde{W}^r(\tilde{\beta}^r_k)$ for
large $r$, and therefore $\tilde{\beta}^r_{k+1}\to kw+w$.
\hfill$\blacksquare$

Let $t\in[0,T]$ be a continuity point of $\sigma(\cdot)$ that is not an
integer multiple of $w$.  Then
$\tilde{\sigma}^r(t)\ws\sigma(t)$. Let
$\mathbf{C}_0\subset\mathbf{C}_b$ be the functionals $f$ with compact support
and $f(0)=0$. For $f\in\mathbf{C}_0$, use \eqref{e.skorohodDynamics} to write
\begin{align*}
  \langle f,\tilde{\sigma}^r(t)\rangle
  & = \left\langle f, \tilde{\mathcal{B}}^r(t)\left(\cdot+_0 
        F^{-1}_{\tilde{\mathcal{B}}^r(t)}(t 
          -\tilde{\beta}^r(t))\right)\right\rangle \\
   &   = \left\langle (1_{(0,\infty)}f)\left(\cdot - 
        F^{-1}_{\tilde{\mathcal{B}}^r(t)}(t-\tilde{\beta}^r(t))\right),
        \tilde{\mathcal{B}}^r(t)\right\rangle.
\end{align*}
Letting $k=\lfloor t/w\rfloor$, Lemma \ref{l.betas} implies that for
all sufficiently large $r$, $t\in(\tilde{\beta}^r_k,\tilde{\beta}^r_{k+1})$.
So for large $r$, $\tilde{\beta}^r(t)=\tilde{\beta}^r_k$ and,
$\tilde{\mathcal{B}}^r(t)=\tilde{\mathcal{E}}^r(\tilde{\beta}^r_k) -
\tilde{\mathcal{E}}^r(\tilde{\beta}^r_{k-1})$ if $t>w$ or
$\tilde{\mathcal{B}}^r(t) = \tilde{\sigma}^r(0)$ if $t<w$. In the first case,
for large $r$,
\[
  \langle f,\tilde{\sigma}^r(t)\rangle
  = \left\langle f\left(\cdot - 
  F^{-1}_{\tilde{\mathcal{E}}^r(\tilde{\beta}^r_k) 
- \tilde{\mathcal{E}}^r(\tilde{\beta}^r_{k-1})}(t-\tilde{\beta}^r_{k-1})\right),
\tilde{\mathcal{E}}^r(\tilde{\beta}^r_k)
-\tilde{\mathcal{E}}^r(\tilde{\beta}^r_{k-1})\right\rangle,
\]
because $1_{(0,\infty)}f=f$. Denote the above integrands by $f^r$. By \eqref{e.skorohodConv},  $\tilde{\mathcal{E}}^r(\cdot)\to \alpha(\cdot)\nu$ uniformly on compact time
intervals in the weak topology and in $L^1$ on $\mathcal{M}_1$. So by Lemma \ref{l.betas}, $\tilde{\mathcal{E}}^r(\tilde{\beta}^r_k)
-\tilde{\mathcal{E}}^r(\tilde{\beta}^r_{k-1})\ws \alpha
w\nu$ and also in $L^1$. Moreover,  $t-\tilde{\beta}^r_k\to\llbracket
t\rrbracket_w$, which is strictly less than $w=\langle\chi,\alpha w\nu\rangle$. 
Thus the assumptions of Lemma \ref{l.Finverse} are satisfied and so
\[ F^{-1}_{\tilde{\mathcal{E}}^r(\tilde{\beta}^r_k) 
- \tilde{\mathcal{E}}^r(\tilde{\beta}^r_{k-1})}(t-\tilde{\beta}^r_{k-1}) 
\to 
F^{-1}_{\alpha w\nu}(\llbracket t \rrbracket_w).
\]
Since $f$ is uniformly continuous, the functions $f^r$ converge uniformly to
$g=f(\cdot - F^{-1}_{\alpha w\nu}(\llbracket t\rrbracket_w)$. Writing
\[
  \langle f,\tilde{\sigma}^r(t)\rangle
  = \left\langle f^r-g,\tilde{\mathcal{E}}^r(\tilde{\beta}^r_k)
-\tilde{\mathcal{E}}^r(\tilde{\beta}^r_{k-1})\right\rangle
+\left\langle g,\tilde{\mathcal{E}}^r(\tilde{\beta}^r_k)
-\tilde{\mathcal{E}}^r(\tilde{\beta}^r_{k-1})\right\rangle,
\]
take $r\to\infty$ on both sides. The first right-hand term is bounded above by 
\[ \|f^r-g\|_\infty\left\langle 1, \tilde{\mathcal{E}}^r(\tilde{\beta}^r_k)
-\tilde{\mathcal{E}}^r(\tilde{\beta}^r_{k-1})\right\rangle
\le  \|f^r-g\|_\infty\tilde{E}^r(T)
\le  \|f^r-g\|_\infty 2\alpha T,
\]
for sufficiently large $r$, which converges to zero. The second right-hand
term converges to $\langle f(\cdot -F^{-1}_{\alpha w\nu}(\llbracket
  t\rrbracket_w),\alpha w\nu\rangle$ since $g$ is bounded and continuous.

\vspace{2ex}
\begin{minipage}[t]{2in}
\footnotesize {\sc Department of Mathematics\\
University of Virginia\\
Charlottesville, VA 22904\\
E-mail:} gromoll@virginia.edu
\end{minipage}
\hspace{5ex}
\begin{minipage}[t]{4in}
\footnotesize {\sc Department of Mathematics\\
SUNY Geneseo\\
Geneseo, NY 14454\\
E-mail:} kochalski@geneseo.edu
\end{minipage}

\end{document}